\documentclass[12pt]{amsart}

\usepackage{amssymb}

\newtheorem{theorem}{Theorem}[section]
\newtheorem{lemma}[theorem]{Lemma}
\newtheorem{proposition}[theorem]{Proposition}
\newtheorem{corollary}[theorem]{Corollary}
\newtheorem{definition}[theorem]{Definition}
\newtheorem{remark}[theorem]{Remark}

\numberwithin{equation}{section}

\begin{document}

\baselineskip=14.5pt

\title[Monodromy group for principal bundle, II]{Monodromy group
for a strongly semistable principal bundle over a curve, II}

\author[I. Biswas]{Indranil Biswas}

\address{School of Mathematics, Tata Institute of Fundamental
Research, Homi Bhabha Road, Mumbai 400005, India}

\email{indranil@math.tifr.res.in}

\author{A. J. Parameswaran}

\address{School of Mathematics, Tata Institute of Fundamental
Research, Homi Bhabha Road, Mumbai 400005, India}

\email{param@math.tifr.res.in}



\date{}

\begin{abstract}

Let $X$ be a geometrically
irreducible smooth projective curve defined over a
field $k$. Assume that $X$ has a $k$--rational point; fix
a $k$--rational point $x\in X$. From these data
we construct an affine group scheme ${\mathcal G}_X$
defined over the field $k$ as well as
a principal ${\mathcal G}_X$--bundle $E_{{\mathcal G}_X}$
over the curve $X$. The group scheme ${\mathcal G}_X$ is given
by a ${\mathbb Q}$--graded neutral
Tannakian category built out of all strongly
semistable vector bundles over $X$. The principal bundle
$E_{{\mathcal G}_X}$ is tautological.
Let $G$ be a linear algebraic group, defined
over $k$, that does not admit any nontrivial character
which is trivial on the connected component, containing
the identity element, of the reduced
center of $G$. Let $E_G$ be a strongly semistable
principal $G$--bundle over $X$. We associate to $E_G$
a group scheme $M$ defined over $k$,
which we call the monodromy group
scheme of $E_G$, and a principal $M$--bundle $E_M$
over $X$, which we call the monodromy bundle of $E_G$.
The group scheme $M$ is canonically
a quotient of ${\mathcal G}_X$, and $E_M$ is the
extension of structure group of $E_{{\mathcal G}_X}$.
The group scheme $M$ is also canonically
embedded in the fiber ${\rm Ad}(E_G)_{x}$
over $x$ of the adjoint bundle.

\end{abstract}

\maketitle

\section{Introduction}

Let $X$ be a geometrically
irreducible smooth projective curve
defined over a field $k$ such that $X$ admits a
$k$--rational point. Fix a $k$--rational
point $x$ of $X$. From this data
we construct a neutral Tannakian
category ${\mathcal C}_X$ defined over the field
$k$ in the following way. The objects of
${\mathcal C}_X$ are all maps $f$ from the rational numbers to
the strongly semistable vector bundles over $X$ such that
$f(\lambda)\,=\, 0$ for all but finitely many $\lambda\,\in\,
{\mathbb Q}$, and if $f(\lambda)\,\not=\, 0$, then
$$
\mu(f(\lambda))\,:=\,
\frac{\text{degree}(f(\lambda))}{\text{rank}(f(\lambda))}\, =\,
\lambda\, .
$$
For any $f, f'\, \in\, {\mathcal C}_X$, set
$$
\text{Hom}(f\, , f')\, :=\, \prod_{\lambda\in {\mathbb Q}}
H^0(X,\, \text{Hom}(f(\lambda)\, , f'(\lambda)))\, .
$$
We note that for any two vector bundles $V$ and $W$
over $X$, with
$$
\mu(V) \, :=\, \frac{\text{degree}(V)}{\text{rank}(V)}\, =\,
\frac{\text{degree}(W)}{\text{rank}(W)} \, =: \mu(W)\, ,
$$
we have $\mu(V\oplus W)\,= \, \mu(V)$.
If $V$ and $W$ are also strongly semistable, then for
any
$$
\phi\,\in\, H^0(X,\, \text{Hom}(V\, , W))\, ,
$$
either the homomorphism
$\phi$ is injective (respectively, surjective) or
$\text{kernel}(\phi)$ (respectively, $\text{cokernel}(\phi)$)
is a strongly semistable vector bundle
with same $\text{degree}/\text{rank}$ quotient as that of $V$;
the details are in Section \ref{sec3}.
These properties imply that ${\mathcal C}_X$ is an abelian
category.

For any $f, f'\, \in\, {\mathcal C}_X$, define their
tensor product
$$
(f\otimes f') (\lambda)\, :=\, \bigoplus_{z\in{\mathbb Q}}
f(z)\otimes f'(\lambda-z)\, ,
$$
and define $f^*$ by $f^*(\lambda) \, =\, f(-\lambda)^*$.
It can be shown that $f\otimes f'\, , f^*\, \in\,
{\mathcal C}_X$.
The object $f_0\, \in\, {\mathcal C}_X$, defined by
$f_0(\lambda)\, =\, 0$ for $\lambda\, \not=\, 0$ and
$f_0(0)\, =\, {\mathcal O}_X$ (the structure sheaf
of $X$), acts as the identity element
for the tensor product operation on ${\mathcal C}_X$.

Using the $k$--rational point $x\, \in\, X$, we have
a fiber functor on ${\mathcal C}_X$ that sends any object
$f$ to the $k$--vector space
$$
\bigoplus_{z\in{\mathbb Q}} f(z)_{x}\, .
$$

All these operations together define a neutral Tannakian
category over $k$.
Let ${\mathcal G}_X$ denote the affine group scheme defined over
$k$ given by this neutral Tannakian category ${\mathcal C}_X$.

Let $\text{Vect}(X)$ denote the category of vector
bundles over the curve $X$.
We have a covariant functor from ${\mathcal C}_X$
to $\text{Vect}(X)$ defined by
$$
f\, \longmapsto\, \bigoplus_{z\in{\mathbb Q}} f(z)\, .
$$
This functor is compatible with the operations of tensor
product, direct sum and dualization. Therefore,
this functor defines a principal
${\mathcal G}_X$--bundle $E_{{\mathcal G}_X}$
over $X$.

The fundamental group scheme of $X$ constructed in \cite{No1},
\cite{No2} is a quotient of ${\mathcal G}_X$.

Let $G$ be a linear algebraic group defined over
the field $k$ with the property that there is no
nontrivial character of $G$ which is trivial on the
center of $G$.
Let $Z_0(G)$ denote the maximal split torus contained
in the reduced center of $G$.

Take a strongly semistable principal $G$--bundle $E_G$ over $X$.
Given a finite dimensional left
$G$--module $V$, consider the isotypical decomposition
$$
V\, =\, \bigoplus_{\chi\in Z_0(G)^*} V_\chi
$$
of the $Z_0(G)$--module $V$, where $Z_0(G)^*$ is the group
of characters of $Z_0(G)$. Since the actions of $Z_0(G)$
and $G$ on $V$ commute, each $V_\chi$ is a $G$--module.
Let $E_V$ (respectively, $E_{V_\chi}$) be the vector bundle
over $X$ associated to the principal $G$--bundle $E_G$ for
the $G$--module $V$ (respectively, $V_\chi$). It can be shown
that the vector bundle $V_\chi$ is strongly semistable. Also,
we have a homomorphism
$$
\delta_{E_G}\, :\, Z_0(G)^*\, \longrightarrow\, {\mathbb Q}
$$
that sends any character $\chi'$ to
$\frac{\text{degree}(V'_{\chi'})}{\text{rank}(V'_{\chi'})}$
for some $G$--module $V'$; it does not depend on the choice
of $V'$ (Corollary \ref{cor1}).

Therefore, for any finite dimensional left
$G$--module $V$, we have an object $f_{E_G,V}$ of ${\mathcal C}_X$ 
defined by
$$
f_{E_G,V}(\lambda)\, :=\, 
\bigoplus_{\{\chi\in Z_0(G)^*\vert \delta_{E_G}(\chi)=\lambda\}}
E_{V_\chi}\, .
$$

We construct a subcategory of the Tannakian category ${\mathcal
C}_X$ by considering all objects of ${\mathcal C}_X$ isomorphic
to some subquotient of some $f_{E_G,V}$, where $V$
runs over all finite dimensional left $G$--modules. This subcategory
gives a quotient group scheme ${\mathcal G}_X$, which we call
the monodromy group scheme of $E_G$. Let $M$ denote the
monodromy group scheme of $E_G$. Let $E_M$ be the principal
$M$--bundle over $X$ obtained by extending the structure
group
of the principal ${\mathcal G}_X$--bundle $E_{{\mathcal G}_X}$.
The details of these constructions are in Section \ref{sec4}.

In \cite{BPS}, the monodromy group scheme and the monodromy bundle
were constructed
under the extra assumptions that the base field is algebraically
closed and the group $G$ is semisimple.

\section{A universal Tannakian category for a
pointed curve}\label{sec3}

Let $k$ be any field. Let $X$ be a geometrically 
irreducible smooth projective curve defined over $k$.

A vector bundle $W$ over $X$ is called \textit{semistable}
if for every subbundle $W'\, \subset\, W$ of positive rank,
the inequality
$\text{degree}(W')/\text{rank}(W')\, \leq\,
\text{degree}(W)\text{rank}(W)$ holds. We recall
that the rational number $\frac{\text{degree}(W)}{\text{rank}(W)}$
is called the \textit{slope} of $W$, and it is denoted by $\mu (W)$.

\begin{proposition}[\cite{La}, \cite{HN}]\label{ba.ch.}
Let $\ell$ be a field extension of $k$. A vector bundle $W$ over $X$
is semistable if and only if the base change $W\bigotimes_k \ell$
over $X\times_k\ell$ is semistable.
\end{proposition}

The above proposition is proved in \cite{La} under the assumption
that $k$ is infinite (see \cite[page 97, Proposition 3]{La}),
and it is proved in \cite{HN} under the assumption
that $k$ is perfect (see \cite[page 222]{HN}). We note that if
$W$ is not semistable, then it is immediate that $W\bigotimes_k \ell$
is not semistable.

Consider the diagram
\begin{equation}\label{ray.}
\begin{matrix}
X &\stackrel{\pi}{\longrightarrow}& X_1
&\stackrel{\varphi}{\longrightarrow}& X\\
&&\downarrow && \downarrow\\
&&\text{Spec}(k) & \stackrel{F_k}{\longrightarrow} & \text{Spec}(k)
\end{matrix}
\end{equation}
where $F_k$
is the Frobenius map of $k$ if the characteristic of
the field $k$ is positive, and it is the identity map
when the characteristic is zero, the
square is Cartesian, and $\pi$ is the
relative Frobenius map (see \cite[page 118]{Ra}). The composition
$\varphi\circ \pi$ will be denoted by $F_X$.

A semistable vector bundle $W$ over $X$ is called \textit{strongly
semistable} if the iterated pull back
$$
(\overset{n-\text{times}}{\overbrace{F_X\circ \cdots \circ F_X}})^*W
$$
is semistable for all $n\,\geq\, 1$.

\begin{remark}\label{r.ba.ch.}
{\rm Strongly semistable vector bundle are usually defined under
the assumption that the base field is perfect. 
In view of Proposition \ref{ba.ch.}, the above definition is
compatible with it.}
\end{remark}

Let ${\mathcal C}_X$ denote the space of all maps $f$ from
$\mathbb Q$ to the space of all strongly semistable
vector bundles over $X$ satisfying the following two
conditions:
\begin{itemize}
\item All but finitely many rational numbers are sent
by $f$ to the vector bundle of rank zero.

\item For any $\eta\,\in\, {\mathbb Q}$ with $f(\eta)\,\not=\, 0$,
$$
\mu(f(\eta))\, :=\,
\frac{\text{degree}(f(\eta))}{\text{rank}(f(\eta))}\, =\, \eta\, .
$$
\end{itemize}

Therefore, ${\mathcal C}_X$ consists of finite collections of
strongly semistable vector bundles
over $X$ of distinct slopes. In other
words, any element of ${\mathcal C}_X$ is of the form
\begin{equation}\label{eq1}
\underline{V}\, :=\, (V_{\lambda_1}\, ,\cdots \, , V_{\lambda_n})\, ,
\end{equation}
where $\lambda_1\, <\, \cdots\, <\, \lambda_n$ are finitely
many (possibly empty)
rational numbers, and for each $i\, \in\, [1\, ,n]$,
$V_i$ is a strongly semistable vector bundle over $X$ with
$$
\mu (V_i) \, =\, \lambda_i\, .
$$

For any $f,f'\,\in\, {\mathcal C}_X$, define $f\bigoplus f'
\,\in\, {\mathcal C}_X$ to be the function that sends
any $\lambda\,\in\, {\mathbb Q}$ to the direct sum of vector
bundles $f(\lambda)\bigoplus f'(\lambda)$. Since both
$f(\lambda)$ and $f'(\lambda)$ are strongly semistable with
$$
\mu (f(\lambda))\, =\, \mu (f'(\lambda))\, =\,
\lambda
$$
provided $f(\lambda)\, \not=\, 0\, \not=\, f'(\lambda)$,
the vector bundle $f(\lambda)\bigoplus f'(\lambda)$ is also
strongly semistable with
$$
\mu(f(\lambda)\oplus f'(\lambda))\,=\, \lambda
$$
provided $f(\lambda)\bigoplus f'(\lambda)\, \not=\, 0$.

If $V$ and $W$ are two strongly semistable vector bundles over
$X$, then the vector bundle $V\bigotimes W$ is also strongly
semistable; this follows from Remark \ref{r.ba.ch.} and
\cite[page 288, Theorem 3.23]{RR}
(reproduced in Theorem \ref{rr.t}). We further have
\begin{equation}\label{eq2}
\mu(V\bigotimes W)\,=\, \mu(V) + \mu(W)\, .
\end{equation}
This enables us to define the tensor product operation
on ${\mathcal C}_X$ in the following way.

For any $f,f'\, \in\, {\mathcal C}_X$, define
\begin{equation}\label{n1}
(f\otimes f') (\lambda)\, :=\, \bigoplus_{z\in{\mathbb Q}}
f(z)\otimes f'(\lambda-z)\, .
\end{equation}
If other words, if $f\, =\, \underline{V}$ as in eqn. \eqref{eq1},
then
$$
(f\otimes f') (\lambda)\, =\, (V_{\lambda_1}\bigotimes
f'(\lambda-\lambda_1))\bigoplus (V_{\lambda_2}\bigotimes
f'(\lambda-\lambda_2)) \bigoplus\cdots \bigoplus
(V_{\lambda_n}\bigotimes f'(\lambda-\lambda_n))\, .
$$
Since both $f(c)$ and $f'(c)$ are zero except for finitely
many $c$, eqn. \eqref{n1} is a finite direct sum. Therefore, using
eqn. \eqref{eq2} we conclude that $f\bigotimes f'\, \in\, {\mathcal C}_X$.

For any $f\, \in\, {\mathcal C}_X$, define the \textit{dual} $f^*$
of $f$
to be the function from $\mathbb Q$ to the strongly semistable vector
bundles over $X$ that sends any $\lambda\, \in\, {\mathbb Q}$ to
the dual vector bundle $f(-\lambda)^*$. Clearly, we have $f^*\, \in\,
{\mathcal C}_X$.

For any $f,f'\, \in\, {\mathcal C}_X$, a \textit{homomorphism} from
$f$ to $f'$ is a function
$$
\gamma\, :\, {\mathbb Q}\, \longrightarrow\, 
\bigoplus_{z\in{\mathbb Q}}H^0(X, \, \text{Hom}(f(z)\, ,f'(z)))
$$
such that $\gamma(z) \, \in\, H^0(X, \, \text{Hom}(f(z)\, ,f'(z)))$
for all $z\,\in\,{\mathbb Q}$. By $\text{Hom}(f\, ,f')$ we will
denote the set of all homomorphisms from $f$ to $f'$. So
$$
\text{Hom}(f\, , f')\, :=\, \prod_{\lambda\in {\mathbb Q}}
H^0(X,\, \text{Hom}(f(\lambda)\, , f'(\lambda)))\, .
$$

A homomorphism $\gamma$ from $f$ to $f'$ will be called an
\textit{isomorphism} if
$$
\gamma(z)\, :\, f(z)\, \longrightarrow\, f'(z)
$$
is an isomorphism for all $z\, \in\, {\mathbb Q}$.

Let $\phi\, :\, V\, \longrightarrow\, W$ be a homomorphism
between strongly semistable vector bundles $V$ and $W$ over $X$
with 
\begin{equation}\label{eq.00}
\mu (V)\, =\, \mu(W) \, .
\end{equation}
Then it can be shown that
either $\phi$ is injective, or $\text{kernel}(\phi)$
is a strongly semistable vector bundle over $X$ with
$$
\mu(\text{kernel}(\phi))\, =\, \mu (V)\, .
$$
Indeed, if $\phi$ is nonzero with $\text{kernel}(\phi)$ a nonzero
subsheaf of $V$, then consider $V/\text{kernel}(\phi)$. We note
that $V/\text{kernel}(\phi)$ is a quotient of $V$ as
well as a subsheaf
of $W$. Therefore, as $V$ and $W$ are semistable, we have
\begin{equation}\label{eq.001}
\mu(V) \, \leq\, \mu (V/\text{kernel}(\phi)) \,\leq\,
\mu(W)\, .
\end{equation}
Now from eqn. \eqref{eq.00} it follows that both the
inequalities in eqn. \eqref{eq.001} are
equalities. Consequently,
$\mu(V) \,=\, \mu (\text{kernel}(\phi))$. It also
follows that $V/\text{kernel}(\phi)$ is torsionfree, because
the inverse image in $V$, of the torsion part of
$V/\text{kernel}(\phi)$, has slope strictly greater than
$\mu ((\text{kernel}(\phi))$ if  $V/\text{kernel}(\phi)$
has torsion. Note that as $V$ is semistable, it does not
have any subsheaf with slope larger than $\mu(V)$. Since
$V/\text{kernel}(\phi)$ is torsionfree, we conclude that
$\text{kernel}(\phi)$ is a subbundle of $V$.

Similarly, either $\phi$ is surjective or
$\text{cokernel}(\phi)\, :=\, W/\phi(V)$ is a
strongly semistable vector bundle over $X$ with
$$
\mu (\text{cokernel}(\phi))\, =\, \mu(W)\, .
$$
(Replace $\phi$ by its dual $\phi^*$ in the above argument.)
This enables us to define the kernel and the cokernel of any
homomorphism between two objects in ${\mathcal C}_X$.

Take any $f,f'\, \in\, {\mathcal C}_X$ and any
$\gamma\, \in\, \text{Hom}(f\, ,f')$. Consider the function from
$\mathbb Q$ to the space of all strongly semistable vector bundles over
$X$ that sends any $\lambda\, \in\, {\mathbb Q}$ to the kernel
of the homomorphism
$$
\gamma (\lambda)\, :\, f(\lambda)\, \longrightarrow\,
f'(\lambda)\, .
$$
This function defines an object of ${\mathcal C}_X$, which
we will call the \textit{kernel} of
$\gamma$, and it will be denoted by $\text{kernel}(\gamma)$.
Similarly, consider the function from $\mathbb Q$ to the
space of all strongly semistable vector bundles over
$X$ that sends any $\lambda\, \in\, {\mathbb Q}$ to the cokernel
of the homomorphism
$$
\gamma (\lambda)\, :\, f(\lambda)\, \longrightarrow\,
f'(\lambda)\, .
$$
The object of ${\mathcal C}_X$ defined by this function will be
denoted by $\text{cokernel}(\gamma)$, and it will be called
the \textit{cokernel} of $\gamma$.

\begin{remark}\label{rem.c}
{\rm The abelian category
$$
{\mathcal C}_X\, =\, \bigoplus_{\lambda\in{\mathbb Q}
} {\mathcal C}^\lambda_X\, ,
$$
where ${\mathcal C}^\lambda_X$ is the abelian category of
strongly semistable vector bundles over $X$ of slope $\lambda$.
The index $\mathbb Q$ in the above direct sum acts as
a weight that guides the tensor product and the dualization
operations.}
\end{remark}

Take any closed point $x$ of $X$. The residue field will be denoted
by $k(x)$. For any vector
bundle $V$ over $X$, the fiber of $V$ over $x$, which is
a $k(x)$--vector space, will be denoted by $V_{x}$.
The
category of finite dimensional vector spaces over the field $k(x)$
will be denoted by $k(x)$--mod.

We have a functor
$$
\omega\, :\, {\mathcal C}_X\, \longrightarrow\, k(x)\mbox{--}{\rm mod}
$$
defined by
$$
f\, \longmapsto\, \bigoplus_{z\in{\mathbb Q}}
f(z)_{x}\, .
$$
In other words, for any
$$
\underline{V}\, :=\, (V_{\lambda_1}\, ,\cdots \, , V_{\lambda_n})
\, \in\, {\mathcal C}_X
$$
as in eqn. \eqref{eq1}, we have $\omega(\underline{V})\,=\,
(V_{\lambda_1})_{x}\bigoplus \cdots \bigoplus
(V_{\lambda_n})_{x}$.

Let ${\mathcal O}_X$ be the trivial line bundle over $X$ defined
by the structure sheaf of $X$. The object in ${\mathcal C}_X$
defined by ${\mathcal O}_X$ will also be denoted by
${\mathcal O}_X$.

The triple $({\mathcal C}_X\, , {\mathcal O}_X\, ,\omega)$
together form a Tannakian category over $k$;
see \cite{Sa}, \cite{DM} for Tannakian category.
We recall that a Tannakian category over $k$ is a rigid abelian
tensor category ${\mathcal C}'$ such that
\begin{itemize}
\item $\text{End}(1\!\!{\mathbf 1}) \, =\, k$, and

\item there is a field extension $k'$ of $k$ and
a $k$--linear fiber functor from ${\mathcal C}'$ to the
category of vector spaces over $k'$.
\end{itemize}

Henceforth, we will assume that $X$ admits a $k$--rational
point. Fix a $k$--rational point $x$ of $X$.

Since $k(x)\, =\, k$, the above Tannakian category
defined by the triple $({\mathcal C}_X\, ,
{\mathcal O}_X\, ,\omega)$ is a
neutral Tannakian category. Hence they define
a affine group scheme defined over $k$ \cite[page 130, Theorem
2.11]{DM}, \cite[Theorem 1.1]{No2}, \cite[Theorem 1]{Sa}.

Let ${\mathcal G}_X$
denote the group scheme defined over $k$ given by the
neutral Tannakian category
$({\mathcal C}_X\, , {\mathcal O}_X\, ,\omega)$.

We note that the above neutral Tannakian category
is ${\mathbb Q}$--graded (see \cite[page 186]{DM}).
We will now show that there is a tautological principal
${\mathcal G}_X$--bundle over $X$.

Take any object
$$
\underline{V}\, :=\, (V_{\lambda_1}\, ,\cdots \, , V_{\lambda_n})
\,\in\, {\mathcal C}_X
$$
as in eqn. \eqref{eq1}. To it we associate the vector bundle
$(V_{\lambda_1}\bigoplus \cdots \bigoplus V_{\lambda_n})$.
In other words, we have a functor
$$
{\mathcal C}_X\, \longrightarrow\, \text{Vect}(X)\, ,
$$
where $\text{Vect}(X)$ as before is the category of vector
bundles over $X$, defined by
\begin{equation}\label{eqm1}
f\, \longmapsto\, \bigoplus_{z\in{\mathbb Q}}
f(z)\, .
\end{equation}
Using \cite[page 149, Theorem 3.2]{DM}, this functor gives a
principal ${\mathcal
G}_X$--bundle over $X$. This principal ${\mathcal G}_X$--bundle
over $X$ will be denoted by $E_{{\mathcal G}_X}$.

\begin{remark}
{\rm A vector bundle $E$ over $X$ is called \textit{finite} if there
are two distinct polynomials $f_1, f_2\, \in\, {\mathbb Z}[t]$
with nonnegative coefficients such that the vector
bundle $f_1(E)$ is isomorphic
to $f_2(E)$. Any finite vector bundle over $X$ is strongly semistable
of degree zero. A vector bundle $V$ of degree zero over $X$ is called
\textit{essentially finite} if there is a finite vector bundle
$E$ over $X$ and a quotient bundle $Q$ of $E$ of degree zero such
that the vector bundle
$V$ is a subbundle of $Q$; see \cite{No1}, \cite{No2}. Essentially
finite vector bundles over $X$ form a neutral Tannakian category.
and the corresponding group scheme is called the \textit{fundamental
group scheme} of $X$ \cite{No2}. Since any essentially finite
vector bundle over $X$ is strongly semistable, the
fundamental group scheme of $X$ is a quotient of the group scheme
${\mathcal G}_X$ constructed above.}
\end{remark}

\section{Monodromy of principal bundles}\label{sec4}

Let $G$ be a linear algebraic group defined over the
field $k$. Let $Z'_0(G)$ denote the connected component,
containing the identity element, of the reduced center of $G$.
We assume that $G$ satisfies the following condition: there is
no nontrivial character of $G$ which is trivial on $Z'_0(G)$.
Let
\begin{equation}\label{eq3}
{Z}_0(G)\, \subset\, Z'_0(G)
\end{equation}
be the (unique) maximal split torus contained in $Z'_0(G)$.
The above condition on $G$ implies that there is
no nontrivial character of $G$ which is trivial on $Z_0(G)$.

The subgroup $Z_0(G)$ gives a decomposition of any $G$--module,
which we will describe now.

Let $V$ be a finite dimensional left $G$--module. Consider the action
on $V$ of the subgroup $Z_0(G)$ in eqn. \eqref{eq3}. Since
$Z_0(G)$ is a product of copies of ${\mathbb G}_m$, the vector
space $V$ decomposes into a direct sum of one dimensional
subspaces such that each of the one dimensional
subspaces is preserved by the action of $Z_0(G)$. Therefore,
we obtain a finite collection of distinct characters of
$Z_0(G)$, say $\chi_1\, ,\cdots
\, , \chi_m$, such that for each line $\xi \, \subset\, V$
preserved by the action of $Z_0(G)$ on $V$, there is a character
$\chi_i$, for some $i\, \in\, [1\, , m]$, such that
$Z_0(G)$ acts on $\xi$ as scalar multiplications
through the character $\chi_i$.

For any $i\, \in\, [1\, , m]$, let
$$
V_i\, \subset\, V
$$
be the linear subspace on which $Z_0(G)$ acts
as scalar multiplications through the
character $\chi_i$. This subspace $V_i$ is preserved by the
action of $G$ on $V$. Indeed, this follows immediately from
the fact that the actions of $G$ and $Z_0(G)$ on $V$ commute.
Therefore, we have a natural decomposition
\begin{equation}\label{eq4}
V\, =\, V_1\bigoplus \cdots \bigoplus V_m
\end{equation}
of the $G$--module $V$. This is clearly the isotypical
decomposition of the $Z_0(G)$-module $V$.

Let $Z_0(G)^*$ denote the group of all characters of $Z_0(G)$.
We may reformulate the decomposition in eqn. \eqref{eq4} in the
following way. Any finite dimensional left $G$--module has
a natural decomposition
\begin{equation}\label{eq5}
V\, =\, \bigoplus_{\chi\in Z_0(G)^*} V_\chi\, ,
\end{equation}
where $V_\chi\, \subset\, V$ is the subspace on which
$Z_0(G)$ acts as scalar multiplications
through the character $\chi$.

As before, let $X$ be a geometrically irreducible smooth
projective curve defined over the field $k$.

\begin{lemma}\label{lem1}
Take any character $\chi\, \in\, Z_0(G)^*$. Let $V$ and
$W$ be two nonzero finite dimensional left $G$--modules such
that $Z_0(G)$ acts on both $V$ and $W$ as scalar
multiplications through the character
$\chi$. Let $E_G$ be a principal $G$--bundle
over $X$. Let $E_V$ (respectively, $E_W$) be the vector bundle
over $X$ associated to the principal $G$--bundle $E_G$
for the $G$--module $V$ (respectively, $W$). Then
$$
\mu(E_V)\, =\, \mu(E_W)\, .
$$
\end{lemma}

\begin{proof}
The $G$--module structures on $V$ and $W$ together induce a
$G$--module structure on the vector space $\text{Hom}(V,W)
\, =\, V^*\bigotimes W$.
The vector bundle over $X$ associated to $E_G$ for the
$G$--module $\text{Hom}(V,W)$ is evidently identified
with the vector bundle $E_V^*\bigotimes E_W$.

Since $Z_0(G)$ acts on both $V$ and $W$ as scalar
multiplications through the
character $\chi$, the action of $Z_0(G)$ on the
$G$--module $\text{Hom}(V,W)$ is the trivial action.
In particular, the action of $Z_0(G)$ on the one--dimensional
$G$--module $\bigwedge^{\text{top}}\text{Hom}(V,W)$
is the trivial action.

By our assumption on $G$, the group
$G/Z_0(G)$ does not admit any nontrivial character.
Therefore, from the observation that
$Z_0(G)$ acts trivially on the one--dimensional
$G$--module $\bigwedge^{\text{top}}\text{Hom}(V,W)$
we conclude that the action of $G$ on
$\bigwedge^{\text{top}}\text{Hom}(V,W)$ is the trivial
action. This immediately implies that the line bundle
over $X$ associated to the principal $G$--bundle $E_G$
for the $G$--module $\bigwedge^{\text{top}}\text{Hom}(V,W)$
is a trivial line bundle.

The line bundle
over $X$ associated to the principal $G$--bundle $E_G$
for the $G$--module $\bigwedge^{\text{top}}\text{Hom}_k(V,W)$
is clearly identified with
$\bigwedge^{\text{top}}(E_V^*\bigotimes E_W)$. Since
$\bigwedge^{\text{top}}(E_V^*\bigotimes E_W)$ is a trivial line
bundle, we have
$$
\text{degree}(\bigwedge^{\text{top}}(E_V^*\bigotimes E_W))\, =\, 0\, .
$$
Now from the identity
$$
\text{degree}(\bigwedge^{\text{top}}(E_V^*\bigotimes E_W))\, =\,
\text{degree}(E_W)\text{rank}(E_V)-
\text{degree}(E_V)\text{rank}(E_W)\, ,
$$
it follows that $\mu(E_V)\, =\, \mu(E_W)$.
This completes the proof of the lemma.
\end{proof}

The following is a corollary of
\cite[page 139, Proposition 2.21]{DM}.

\begin{corollary}\label{ch.ar.}
All characters of $Z_0(G)$ arise from the indecomposable
representations of $G$.
In other words, for any character $\chi$ of $Z_0(G)$,
there is some nonzero finite dimensional indecomposable
left $G$--module $V$ such that $Z_0(G)$ acts on $V$ as
scalar multiplications through the character $\chi$.
\end{corollary}

If $E$ and $F$ are two vector bundles over $X$, then
$\mu(E\bigotimes F)\, =\, \mu(E)+   \mu(F)$.
Therefore, Lemma \ref{lem1} and Corollary \ref{ch.ar.}
combine together to give the following corollary:

\begin{corollary}\label{cor1}
Fix any principal $G$--bundle $E_G$ over $X$.
Then there is a homomorphism to the additive group
$$
\delta_{E_G}\, :\, Z_0(G)^*\, \longrightarrow\, {\mathbb Q}
$$
that sends any character $\chi$ to
$\frac{{\rm degree}(E_V)}{{\rm rank}(E_V)}$, where
$V$ is a finite dimensional nonzero left $G$--module on which
$Z_0(G)$ acts as scalar multiplications
through the character $\chi$, and $E_V$ is the vector bundle
over $X$ associated to the principal $G$--bundle $E_G$
for the $G$--module $V$.
\end{corollary}

\begin{definition}\label{st.st.}
{\rm Let $G$ be any affine group scheme defined over $k$. A
principal $G$--bundle $E_G$ over a geometrically irreducible
smooth projective curve $X$ will be called \textit{strongly
semistable} if for any indecomposable finite dimensional
left $G$--module $V\, \in\, G\mbox{--}{\rm mod}$, the
vector bundle over $X$ associated to $E_G$ for $V$ is
strongly semistable.}
\end{definition}

See \cite{Ra}, \cite{RR} for the definition of a
(strongly) semistable principal bundle with
a reductive group as the structure group. We will
show that the above definition
coincides with the usual definition when $G$ is reductive.
For that we will need the following theorem.

\begin{theorem}[RR, Theorem 3.23]\label{rr.t}
Let $H$ and $H'$ be reductive linear algebraic groups
defined over $k$ and
$$
\rho\, :\, H\, \longrightarrow\, H'
$$
a homomorphism of algebraic groups such that
$\rho(Z'_0(H))\, \subset\, Z'_0(H')$, where
$Z'_0(H)$ (respectively, $Z'_0(H')$)
is the connected component, containing the identity
element, of the reduced center of $H$ (respectively,
$H'$). Let $E_H$ be
a strongly semistable principal $H$--bundle over $X$.
Then the principal $H'$--bundle $E_{H'}\, :=\, E_H(H')$
over $X$, obtained by extending the structure group
of $E_H$ using $\rho$, is also strongly semistable.
\end{theorem}

Definition \ref{st.st.} is justified by the following
lemma.

\begin{lemma}\label{s.c.}
Let $H$ be a reductive linear algebraic group
and $X$ a geometrically irreducible smooth projective curve
defined over $k$. A principal
$H$--bundle $E_H$ over $X$ is strongly semistable if and only
if for every indecomposable $H$--module $V$, the vector bundle
$E_V\, =\, E_H(V)$ over $X$ associated to the principal
$H$--bundle $E_H$ for $V$ is strongly semistable.
\end{lemma}

\begin{proof}
Let $E_H$ be a strongly semistable principal $H$--bundle
over $X$.
Take any indecomposable $H$--module $V$. Since
$V$ is indecomposable, the group $Z'_0(H)$ acts on $V$
as scalar multiplications, where $Z'_0(H)$ as before
is the connected component, containing the identity element,
of the reduced center of $H$. Therefore, from
Theorem \ref{rr.t} it follows that
the associated vector bundle $E_V$ is strongly semistable.

To prove the converse, let $E_H$ be principal $H$--bundle over $X$
such that for every indecomposable $H$--module $V$, the vector bundle
over $X$ associated to the principal
$H$--bundle $E_H$ for $V$ is strongly semistable.

Let ${\rm ad}(E_H)$ be the adjoint vector bundle for $E_H$.
We recall that ${\rm ad}(E_H)$ is the vector bundle
over $X$ associated to the principal $H$--bundle $E_H$
for the adjoint action of $H$ on its Lie algebra
$\mathfrak h$. To prove
that the principal $H$--bundle $E_H$ is strongly semistable,
it suffices to show that the vector bundle ${\rm ad}(E_H)$ is
strongly semistable. To see this, let
$$
E_P\, \subset\, (F^r_X)^*E_H
$$
be a reduction of structure group to a maximal parabolic
subgroup $P\, \subset\, H$ that violates the semistability
condition, or in other words, we have
\begin{equation}\label{equu}
\text{degree}(((F^r_X)^* \text{ad}(E_H))/\text{ad}(E_P))
\, <\, 0
\end{equation}
(see \cite{Ra}, \cite{RR}). Consider the subbundle
$$
\text{ad}(E_P)\, \subset\, \text{ad}((F^r_X)^*E_H)\, =\,
(F^r_X)^* \text{ad}(E_H)\, .
$$
Since $\text{degree}(\text{ad}(E_H))\, =\, 0$ (as
$\bigwedge^{\rm top} {\mathfrak h}$ is a trivial
$H$--module), from eqn. \eqref{equu}
it follows immediately that the subbundle
$\text{ad}(E_P)$ contradicts the semistability condition
of the vector bundle $(F^r_X)^* \text{ad}(E_H)$.

Consider the Lie algebra $\mathfrak h$ as a $H$--module
using the adjoint action. Note that $Z'_0(H)$ acts trivially
on $\mathfrak h$. Let
\begin{equation}\label{deq}
{\mathfrak h}\, =\, \bigoplus_{i=1}^n V_i
\end{equation}
be a decomposition of ${\mathfrak h}$ into a direct sum of
indecomposable $H$--modules. Let $E_{V_i}$ denote the
vector bundle over $X$ associated to the principal
$H$--bundle $E_H$ for the $H$--module $V_i$. From eqn.
\eqref{deq} we have
$$
\text{ad}(E_H)\, =\,\bigoplus_{i=1}^n E_{V_i}\, .
$$
{}From the given condition on $E_H$ we know that
$E_{V_i}$ is strongly semistable for all $i\, \in\,
[1\, ,n]$.

As $Z'_0(H)$ acts trivially on $\mathfrak h$, and
$H/Z_0(H)$ being semisimple does not have any nontrivial
character, the induced action of $H$ on the line
$\bigwedge^{\rm top} V_i$ is the trivial action. Therefore,
$\bigwedge^{\rm top}E_{V_i}$ is a trivial line bundle.
In particular,
$$
\text{degree}(E_{V_i}) \, =\, 0
$$
for all $i\, \in\, [1\, ,n]$.

Since $\text{ad}(E_H)$ is a direct sum of strongly
semistable vector bundles of degree zero, we conclude
that the vector bundle $\text{ad}(E_H)$ is strongly
semistable. This completes the proof of the lemma.
\end{proof}

Let $E_G$ be a strongly semistable principal $G$--bundle
over $X$. To each $G$--module we will associate an object
of the neutral
Tannakian category ${\mathcal C}_X$ that we constructed in
Section \ref{sec3}.

Let $V$ be a finite dimensional left $G$--module. First consider
the natural decomposition into a direct sum of $G$--modules
$$
V\, =\, \bigoplus_{\chi\in Z_0(G)^*} V_\chi
$$
constructed in eqn. \eqref{eq5}. Let $E_{V_\chi}$ be the
vector bundle over $X$ associated to the principal $G$--bundle
$E_G$ for the above $G$--module $V_\chi$.

\begin{lemma}\label{lem2}
The vector bundle $E_{V_\chi}$ is strongly semistable, and
if $V_\chi\, \not=\, 0$, then
$$
\mu (E_{V_\chi})\,=\, \delta_{E_G}(\chi)\, ,
$$
where $\delta_{E_G}$ is the homomorphism constructed in
Corollary \ref{cor1}.
\end{lemma}

\begin{proof}
Expressing $V_\chi$ as a direct sum of indecomposable
$G$--modules, and using the fact that a direct sum
of strongly semistable vector bundles of same slope
remains strongly semistable, we conclude
that the associated vector bundle
$E_{V_\chi}$ is strongly semistable
(see Definition \ref{st.st.}).

{}From the definition of $\delta_{E_G}$ it follows immediately
that $\mu(E_{V_\chi}) \,=\, \delta_{E_G}(\chi)$
if $V_{\chi}\, \not=\, 0$. This completes the proof of the lemma.
\end{proof}

Lemma \ref{lem2} has the following corollary:

\begin{corollary}\label{cor2}
For any $\lambda\, \in\, {\mathbb Q}$ and any
$V\,\in\, G\mbox{--}{\rm mod}$, the direct sum
$$
E^\lambda_G(V)\, :=\,
\bigoplus_{\{\chi\in Z_0(G)^*\vert \delta_{E_G}(\chi)=\lambda\}}
E_{V_\chi}
$$
is either zero, or it is a strongly semistable vector bundle
with
$$
\mu(E^\lambda_G(V))\,=\, \lambda\, .
$$
\end{corollary}

Finally, consider the function $f_{E_G,V}$ from $\mathbb Q$ to the
space of all vector bundles over $X$ defined 
by
\begin{equation}\label{eq6}
f_{E_G,V}(\lambda)\, :=\, E^\lambda_G(V)\, ,
\end{equation}
where $V\,\in\, G\mbox{--}{\rm mod}$, and
$E^\lambda_G(V)$ is defined in Corollary \ref{cor2}. From
Corollary \ref{cor2} it follows immediately that
this function $f_{E_G,V}$ is an object of
the category ${\mathcal C}_X$ constructed in Section \ref{sec3}.
Therefore, to each object of $G$--mod we have associated
an object of ${\mathcal C}_X$.

Our aim is to construct a Tannakian category out of the
strongly semistable principal $G$--bundle $E_G$. Before that
we will introduce some definitions.

Take any object $f$ in the category ${\mathcal C}_X$.
A \textit{sub--object}
of $f$ is an object $f'$ in ${\mathcal C}_X$ such that
for each $\lambda \, \in\, {\mathbb Q}$, the vector bundle
$f'(\lambda)$ is a subbundle of the vector bundle
$f(\lambda)$. If $f'$ is a sub--object of $f$, then the
object of ${\mathcal C}_X$ that sends any
$\lambda \, \in\, {\mathbb Q}$
to the quotient vector bundle $f(\lambda)/f'(\lambda)$ will
be called a \textit{quotient--object} of $f$.

Let $E$ be a strongly semistable vector bundle over $X$
and $E'$ a nonzero proper subbundle of $E$ with
$$
\mu(E') \, =\, \mu(E)\, .
$$
Then $E'$ is strongly semistable, and furthermore, the quotient
vector bundle $E/E'$ is also strongly semistable with
$\mu(E/E') \, =\,\mu(E)$
if $E/E'\, \not=\, 0$.
Therefore, for an object $f$ of the category ${\mathcal C}_X$,
any quotient--object of $f$ also lie in ${\mathcal C}_X$. Also,
given a subbundle $V'_\lambda\, \subset\, f(\lambda)$ for
each $\lambda\, \in\, {\mathbb Q}$, to check that the function
$$
\lambda\, \longmapsto\, V'_\lambda
$$
is a sub--object of $f$, all we need to check that
$$
\mu(V'_\lambda)\, =\,\lambda
$$
for all $\lambda\, \in\, {\mathbb Q}$ with
$V'_\lambda\, \not=\, 0$.

For any object $f$ of the category ${\mathcal C}_X$, a
\textit{sub-quotient} of $f$ is a sub--object of some
quotient--object of $f$.

Let ${\mathcal C}_{E_G}$ denote the subcategory of ${\mathcal C}_X$
defined by all objects $f$ of ${\mathcal C}_X$ such that there
exists some $V\, \in\, G\mbox{--}{\rm mod}$
with the property that $f$ is isomorphic
to a sub-quotient of $f_{E_G,V}$, where $f_{E_G,V}$ is the object
of ${\mathcal C}_X$ constructed from $V$
in eqn. \eqref{eq6}. The morphisms remain unchanged. In other
words, for any two objects $f$ and $f'$ in ${\mathcal C}_{E_G}$,
the morphisms from $f$ to $f'$ are the morphisms from $f$ to $f'$
considered as objects of ${\mathcal C}_X$.

It is straight--forward to check that
${\mathcal C}_{E_G}$ is a neutral Tannakian subcategory of ${\mathcal
C}_X$. Therefore, the neutral Tannakian
category ${\mathcal C}_{E_G}$ gives an
affine group scheme defined over $k$.

\begin{definition}\label{def01}
{\rm The affine group scheme defined over $k$ given by the
neutral Tannakian
category ${\mathcal C}_{E_G}$ will be called the \textit{monodromy
group scheme} of $E_G$. The monodromy group scheme of $E_G$ will
be denoted by $M$.}
\end{definition}

Since ${\mathcal C}_{E_G}$ is a Tannakian subcategory of
${\mathcal C}_X$, the monodromy group scheme $M$ is a
quotient of the group scheme
${\mathcal G}_X$ constructed in Section \ref{sec3}
(see \cite[Proposition 2.21]{DM}). Let
\begin{equation}\label{mo.-de.}
\phi_{E_G} \, :\, {\mathcal G}_X\, \longrightarrow\, M
\end{equation}
be the quotient map.

Just as we have the tautological ${\mathcal G}_X$--bundle
$E_{{\mathcal G}_X}$ (see eqn. \eqref{eqm1}), there is a tautological
principal $M$--bundle over $X$.

\begin{definition}\label{def02}
{\rm Let $E_M$ denote the tautological principal $M$--bundle
over $X$. This principal $M$--bundle
$E_M$ will be called the \textit{monodromy bundle} for $E_G$.}
\end{definition}

The principal $M$--bundle $E_M$ is evidently the one
obtained by extending the structure group
of the principal ${\mathcal G}_X$--bundle $E_{{\mathcal G}_X}$
using the homomorphism $\phi_{E_G}$ in eqn. \eqref{mo.-de.}.

We will next show that there is a tautological embedding of
the monodromy group scheme $M$ into
the fiber, over the fixed $k$--rational point $x \,\in\,X$,
of the adjoint bundle for $E_G$.

Let $\text{Ad}(E_G)$ be the adjoint bundle for the
principal $G$--bundle $E_G$
over $X$. Let $\text{Ad}(E_G)_x$ be the fiber
of $\text{Ad}(E_G)$ over the fixed $k$--rational point $x$ of $X$.

If $\omega$ is the fiber functor for the principal
$G$--bundle $E_G$ over $X$, then the group
$\text{Ad}(E_G)_x$ defined over $k$
represents the functor $\underline{\text{Aut}}^{\otimes}(\omega)$.
Using Theorem 2.11 in \cite[page 130]{DM}, we get a natural
homomorphism from the group scheme ${\mathcal G}_X$
(constructed in Section \ref{sec3}) to
$\text{Ad}(E_G)_x$. Let
\begin{equation}\label{p.e.g}
\Phi(E_G)\, :\, {\mathcal G}_X \, \longrightarrow\,
\text{Ad}(E_G)_x
\end{equation}
be this natural homomorphism.

It is easy to see that the image of the homomorphism
$\Phi(E_G)$ in eqn. \eqref{p.e.g} coincides with the
monodromy group scheme $M$ in Definition \ref{def01}.

Therefore, we have the following proposition:

\begin{proposition}\label{prop.id.}
The monodromy group scheme $M$ for $E_G$ (introduced in
Definition \ref{def01}) is identified with the image
of the homomorphism $\Phi(E_G)$ in eqn. \eqref{p.e.g}. In other
words, the kernel of the homomorphism $\Phi(E_G)$ coincides
with the kernel of the homomorphism $\phi_{E_G}$ in
eqn. \eqref{mo.-de.}.

There is a natural inclusion $M\, \hookrightarrow\,
{\rm Ad}(E_G)_{x}$ obtained from the fact that the
quotients of ${\mathcal G}_X$ for the two homomorphisms
$\Phi(E_G)$ and $\phi_{E_G}$ coincide.
\end{proposition}

We will now investigate the behavior of the monodromy
group and the monodromy bundle under the extensions of
structure group.

Let
\begin{equation}\label{rh.}
\rho\,:\, G\, \longrightarrow\, G_1
\end{equation}
be an algebraic homomorphism between linear algebraic groups
defined over $k$. Let $Z_0(G_1)$ denote the (unique) maximal
split torus contained in the reduced center of $G_1$. We assume
the following:
\begin{itemize}
\item The group $G_1$ does not admit any nontrivial character
which is trivial on $Z_0(G_1)$.

\item The homomorphism $\rho$ in eqn. \eqref{rh.} satisfies
the condition
\begin{equation}\label{inc.}
\rho(Z_0(G))\, \subset\, Z_0(G_1)\, .
\end{equation}
\end{itemize}

\begin{lemma}\label{rr.g.l.}
Let $E_G$ be a strongly semistable principal
$G$--bundle over $X$. Then the principal $G_1$--bundle
$E_{G_1} \, :=\, E_G(G_1)$, obtained by extending
the structure group of $E_G$ using $\rho$ (defined
in eqn. \eqref{rh.}) is also strongly semistable.
\end{lemma}

\begin{proof}
Take an indecomposable $G_1$--module $V$. So
$Z_0(G_1)$ acts on $V$ as scalar multiplications
through a character. Therefore, using eqn. \eqref{inc.}
we know that $Z_0(G)$ acts on $V$ as scalar
multiplications through a character. Since $E_G$
is strongly semistable, Corollary \ref{cor2} says that
the associated vector bundle
$E_G(V) \,=\, E_{G_1}(V)$ is
strongly semistable. Hence $E_{G_1}$ is strongly
semistable (see Definition \ref{st.st.}).
\end{proof}

Let $E_G$ be a strongly semistable principal
$G$--bundle over $X$. Hence by Lemma \ref{rr.g.l.},
the principal $G_1$--bundle $E_{G_1}$, obtained by
extending the structure group of the principal $G$--bundle
$E_G$ using the homomorphism $\rho$, is also
strongly semistable.
The homomorphism $\rho$ in eqn. \eqref{rh.}
induces a homomorphism of group schemes
\begin{equation}\label{eq4.}
\widetilde{\rho}\, :\, \text{Ad}(E_G)\, \longrightarrow\,
\text{Ad}(E_{G_1})
\end{equation}
over $X$.

\begin{lemma}\label{le1}
The monodromy group scheme $M_1\, \subset\,
{\rm Ad}(E_{G_1})_x$ for $E_{G_1}$
is the image $\widetilde{\rho}(x)(M)$, where
$\widetilde{\rho}(x)$ is the homomorphism in
eqn. \eqref{eq4.} restricted to the $k$--rational
point $x$ of $X$, and $M\, \subset\, {\rm Ad}
(E_G)_{x}$ is the monodromy group scheme
of $E_G$ (see Proposition \ref{prop.id.}). Furthermore,
the monodromy bundle $E_{M_1}$ for $E_{G_1}$ is the
extension of structure group of the monodromy bundle
$E_{M}$ for $E_{G}$ by the homomorphism
$M\, \longrightarrow\, M_1$ obtained by
restricting $\widetilde{\rho}(x)$.
\end{lemma}

\begin{proof}
Take any $G_1$--module $V\, \in\, G_1\mbox{--}{\rm mod}$.
The $G$--module given by $V$ using the homomorphism $\rho$
in eqn. \eqref{rh.} 
will also be denoted by $V$. Consider the isotypical
decomposition of the $Z_0(G)$--module $V$, and also
consider the isotypical
decomposition of the $Z_0(G_1)$--module $V$.
The second
decomposition is finer in the following sense. From
eqn. \eqref{inc.} we get a homomorphism
$$
\rho^*\, :\, Z_0(G_1)^*\, \longrightarrow\, Z_0(G)^*
$$
of character groups.
For any $\chi\, \in\, Z_0(G)^*$, the component of $V$
on which $Z_0(G)$ acts as scalar multiplications through
$\chi$ is the direct sum
$$
\bigoplus_{\chi'\in\, (\rho^*)^{-1}(\chi)} V_{\chi'}\, ,
$$
where $ V_{\chi'}\, \subset\, V$ is the subspace on which
$Z_0(G_1)$ acts as scalar multiplications through
$\chi'$. 

Using this observation it follows that the neutral Tannakian 
category ${\mathcal C}_{E_{G_1}}$ for the principal
$G_1$--bundle $E_{G_1}$ is a subcategory of the neutral
Tannakian category ${\mathcal C}_{E_G}$ for $E_G$.
Now the lemma follows from the
constructions of the monodromy group scheme and the
monodromy bundle and the criterion for surjectivity
in \cite[page 139, Proposition 2.21(a)]{DM}.
\end{proof}

\section{Some properties of the
monodromy group scheme}\label{se5}

As before, let $G$ be a linear algebraic group
defined over $k$ which does not admit any nontrivial
character trivial on $Z_0(G)$.
Take a principal $G$--bundle $E_G$ over the
curve $X$.

For a subgroup scheme $H\, \subset\, G$, let
$(H\bigcap Z_0(G))_{\rm red}$ denote the reduced intersection;
let $(H\bigcap Z_0(G))_0$ denote the (unique) maximal
split torus contained in the abelian
group $(H\bigcap Z_0(G))_{\rm red}$.

\begin{definition}\label{de.bal.}
{\rm A reduction of structure group
\begin{equation}\label{red.}
E_H\, \subset\, E_G
\end{equation}
of $E_G$ to a subgroup scheme $H\, \subset\, G$ will
be called \textit{balanced} if for every character
$\chi$ of $H$ trivial on $(H\bigcap Z_0(G))_0
\, \subset\, H$ (see the above definition), we have
$$
\text{degree}(E_H(\chi))\, =\, 0\, ,
$$
where $E_H(\chi)$ is the line bundle over $X$ associated
to the principal $H$--bundle $E_H$ for the character
$\chi$.}
\end{definition}

\begin{remark}\label{rem2}
{\rm Since any character of $(H\bigcap Z_0(G))_{\rm red}/
(H\bigcap Z_0(G))_0$ is of finite order,
if a character $\chi$ of $H$ is trivial on
$(H\bigcap Z_0(G))_0$, then there is a positive integer $n$
such that the character $\chi^n$ of $H$ is trivial on
$(H\bigcap Z_0(G))_{\rm red}$. Therefore, a reduction
$E_H\, \subset\, E_G$ as in Definition \ref{de.bal.} is balanced
if and only if for every character $\chi$ of $H$ trivial on
$(H\bigcap Z_0(G))_{\rm red}$ we have
$$
\text{degree}(E_H(\chi))\, =\, 0\, .
$$

Since the quotient $(H\bigcap Z_0(G))/(H\bigcap Z_0(G))_{\rm red}$
is a finite group scheme, if a character $\chi$ of $H$ is trivial on
$(H\bigcap Z_0(G))_{\rm red}$, then there is a positive integer $n$
such that the character $\chi^n$ of $H$ is trivial on $H\bigcap Z_0(G)$.
Therefore, a reduction
$E_H\, \subset\, E_G$ as in Definition \ref{de.bal.} is balanced
if and only if for every character $\chi$ of $H$ trivial on
$H\bigcap Z_0(G)$ we have
$$
\text{degree}(E_H(\chi))\, =\, 0\, .
$$}
\end{remark}

\begin{proposition}\label{com.}
Let $E_G$ be a strongly semistable principal $G$--bundle
over $X$ and $E_H\, \subset\, E_G$ a
balanced reduction of structure group of $E_G$
to a subgroup scheme $H\, \subset\, G$.
Then the principal $H$--bundle
$E_H$ over $X$ is strongly semistable.
\end{proposition}

\begin{proof}
Take any indecomposable $H$--module $W$. Let $E_W\, =\,
E_H(W)$ be the vector bundle over $X$ associated to
the principal $H$--bundle $E_H$ for the $H$--module
$W$. We need to show that $E_W$ is strongly semistable.
For that it suffices to show that the vector bundle
${\rm End}(E_W)$ is strongly semistable. Indeed, if
a subbundle $F\, \subset\, (F^j_X)^*E_W$ contradicts the
semistability condition of the vector
bundle $(F^j_X)^*E_W$, then the subbundle
$$
((F^j_X)^*E_W/F)^*\bigotimes (F^j_X)^*E_W\, \subset\,
((F^j_X)^*E_W)^*\bigotimes (F^j_X)^*E_W \, =\, 
(F^j_X)^*{\rm End}(E_W)
$$
contradicts the semistability condition of
the vector bundle $(F^j_X)^*{\rm End}(E_W)$.

Let
\begin{equation}\label{rhoc}
\rho\, :\, H\, \longrightarrow\, \text{GL}(\text{End}(W))
\end{equation}
be the homomorphism given by the action of $H$
on $\text{End}(W)$ induced by the action of $H$
on $W$. We note that ${\rm End}(E_W)$
is the vector bundle associated to
the principal $H$--bundle $E_H$ for the $H$--module
$\text{End}(W)$.

Since the $H$--module $W$ is indecomposable, the group
scheme
$H\bigcap Z_0(G)$ acts on $W$ as scalar multiplications
through a character of $H\bigcap Z_0(G)$.
Therefore, $H\bigcap Z_0(G)$ acts trivially on
$\text{End}(W)$. In other words, $\text{End}(W)$ is an
$H/(H\bigcap Z_0(G))$--module.

As $H$ is a subgroup scheme of $G$, we have
$$
H/(H\bigcap Z_0(G))\, \subset\, G/Z_0(G)\, .
$$
Therefore, there
is a $G/Z_0(G)$--module $V$ such that the
$H/(H\bigcap Z_0(G))$--module $\text{End}(W)$
is a subquotient of the $H/(H\bigcap Z_0(G))$--module
$V$ (see \cite[page 139, Proposition 2.21(b)]{DM}). Let
\begin{equation}\label{q}
V\, \twoheadrightarrow\, Q
\end{equation}
be a quotient of the $H/(H\bigcap Z_0(G))$--module
$V$ such that $\text{End}(W)$ is a submodule of
the $H/(H\bigcap Z_0(G))$--module $Q$.

Let
\begin{equation}\label{d.s.d}
V\, =\, \bigoplus_{i=1}^\ell V_i
\end{equation}
be a decomposition of the $G/Z_0(G)$--module $V$
into a direct sum of indecomposable $G/Z_0(G)$--modules.
For any $i\, \in\, [1\, ,\ell]$, let
$E_{V_i}$ be the vector bundle over $X$ associated to
the principal $G$--bundle $E_G$ for the $G$--module
$V_i$ (the $G/Z_0(G)$--module $V_i$ is considered
as a $G$--module using the quotient map to $G/Z_0(G)$).

Since $E_G$ is strongly semistable, and the $G$--module
$V_i$ is indecomposable, the associated vector bundle
$E_{V_i}$ is strongly semistable. As $G/Z_0(G)$ does
not admit any nontrivial characters, the induced
action of $G$ on the line $\bigwedge^{\rm top}V_i$
is the trivial action. Therefore, the
associated line bundle
$\bigwedge^{\rm top}E_{V_i}$ is trivializable. In
particular, we have
$$
\text{degree}(E_{V_i})\, =\, 0\, .
$$
Since each $E_{V_i}$ is strongly semistable of degree
zero, the vector bundle $E_V$ is also
strongly semistable of degree zero.

Let $E_Q$ denote the vector bundle over $X$
associated to the principal $H$--bundle $E_H$ for the
$H$--module $Q$ in eqn. \eqref{q} (the
$H/(H\bigcap Z_0(G))$--module $Q$ is considered
as an $H$--module using the quotient map
to $H/(H\bigcap Z_0(G))$). Since $H\bigcap Z_0(G)$
acts trivially on the line $\bigwedge^{\rm top}Q$,
and $E_H\, \subset\, E_G$ is a balanced reduction
of structure group, we have
$$
\text{degree}(E_Q)\, =\, \text{degree}(E_H
(\bigwedge^{\rm top}Q))\, =\, 0\, ,
$$
where $E_H(\bigwedge^{\rm top}Q)$ is the line bundle over
$X$ associated to the principal $H$--bundle $E_H$ for the
$H$--module $\bigwedge^{\rm top}Q$.

Since $Q$ is a quotient of the $H$--module $V$,
the vector bundle $E_Q$ is a quotient of $E_V$.
The vector bundle $E_V$ is strongly semistable
of degree zero, and $E_Q$ is a quotient of it
of degree zero. Hence the vector bundle $E_Q$
is also strongly semistable.

We recall that the $H$--module $\text{End}(W)$ is a
submodule of the $H$--module $Q$. Therefore,
the associated vector bundle ${\rm End}(E_W)$ is
a subbundle of $E_Q$. Since $E_Q$ is a strongly semistable
vector bundle of degree zero, and $\text{End}(E_W)$ is a
subbundle of it of degree zero, we conclude that the
vector bundle $\text{End}(E_W)$ is strongly semistable.

We saw earlier that $E_W$ is strongly semistable
if $\text{End}(E_W)$ is so. Therefore, the principal
$H$--bundle $E_H$ is strongly semistable. This
completes the proof of the proposition.
\end{proof}

Let $E_G$ be a strongly semistable principal $G$--bundle
over $X$. In Proposition \ref{prop.id.} we saw that
the monodromy group scheme $M$ (constructed in Definition
\ref{def01}) is canonically embedded in
$\text{Ad}(E_G)_{x}$. For national convenience,
we will denote by $\widetilde{G}$ the group
$\text{Ad}(E_G)_{x}$ defined over
$k$. Let $E_{\widetilde{G}}$ be the principal
$\widetilde{G}$--bundle over $X$ obtained
by extending the structure group of the monodromy bundle
$E_M$ (see Definition \ref{def02}) using the inclusion
of $M$ in $\widetilde{G}$. Therefore,
\begin{equation}\label{red.nv.}
E_M\,\subset\, E_{\widetilde{G}}
\end{equation}
is a reduction
of structure group of $E_{\widetilde{G}}$ to $M$.

Let $Z_0(\widetilde{G})$ denote the unique maximal
split torus contained in the reduced center of
$\widetilde{G}$.

\begin{theorem}\label{thm1}
Assume that the group $\widetilde{G}\,:=\,
\text{Ad}(E_G)_{x}$ does not admit
any nontrivial character which is trivial on
$Z_0(\widetilde{G})$. Then the 
reduction of structure group in eqn. \eqref{red.nv.}
is a balanced reduction of structure
group of $E_{\widetilde{G}}$ to $M$. In particular, the
principal $M$--bundle $E_M$ is strongly semistable.
\end{theorem}

\begin{proof}
We first note that the quotient $M/(M\bigcap Z_0(\widetilde{G}))$
is a subgroup scheme of
$$
G'\, :=\, \widetilde{G}/Z_0(\widetilde{G})\, .
$$
Therefore,
any $M/(M\bigcap Z_0(\widetilde{G}))$--module is a subquotient of some
$G'$--module considered as a $M/(M\bigcap Z_0(G))$--module
(see  \cite[page 139, Proposition 2.21(b)]{DM}).
Let $\chi$ be a character of $M$ which is trivial on
the group scheme $M\bigcap Z_0(\widetilde{G})$. Let $L$ denote
the line bundle over $X$ associated to the
principal $M$--bundle $E_M$ for the character
$\chi$. To prove that $E_M\, \subset\, E_{\widetilde{G}}$ is a
balanced reductive of structure group, it suffices
to show that $\text{degree}(L)\, =\, 0$
(see Remark \ref{rem2}).

The one--dimensional
$M$--module corresponding to $\chi$ will be denoted
by $\xi$. Since $\chi$ is trivial on $M\bigcap Z_0(\widetilde{G})$,
the $M$--module $\xi$ is given by a
$M/(M\bigcap Z_0(\widetilde{G}))$--module. This
$M/(M\bigcap Z_0(\widetilde{G}))$--module will also be denoted by
$\xi$. Let $V$ be a $G'$--module such that the
$M/(M\bigcap Z_0(\widetilde{G}))$--module $\xi$ is a subquotient of
$V$ (we noted earlier that any $M/(M\bigcap
Z_0(\widetilde{G}))$--module is a subquotient of some $G'$--module).

Since $E_G$ is strongly semistable, the $\widetilde{G}$--bundle
$E_{\widetilde{G}}$ is also strongly semistable.
Let $E_V$ denote the vector bundle over $X$ associated
to the principal $\widetilde{G}$--bundle $E_{\widetilde{G}}$
for the $\widetilde{G}$--module
$V$ (since $G'$ is a quotient of $\widetilde{G}$, any $G'$--module
is also a $\widetilde{G}$--module).
As $Z_0(\widetilde{G})$ acts trivially on $V$,
and the principal $\widetilde{G}$--bundle $E_{\widetilde{G}}$
is strongly semistable,
the associated vector bundle $E_V$ is also strongly semistable
(see Lemma \ref{rr.g.l.}).

By our assumption, $G'$ does not admit any nontrivial
character. Therefore, the one--dimensional $\widetilde{G}$--module
$\bigwedge^{\rm top}V$ is a trivial $\widetilde{G}$--module.
Consequently, the line bundle over $X$ associated
to the principal $\widetilde{G}$--bundle $E_{\widetilde{G}}$ for the
one--dimensional $\widetilde{G}$--module
$\bigwedge^{\rm top}V$ is a trivializable.
Hence, we have
$$
\text{degree}(E_V) \, =\, \text{degree}(\bigwedge^{\rm top}E_V) 
\, =\, 0\, .
$$
The earlier defined line bundle $L$ over $X$ is the one associated
to the principal $M$--bundle $E_M$ for the $M$--module
$\xi$. From the definition of the monodromy
bundle $E_M$ it follows immediately that $L$
is an object of the category ${\mathcal C}_X$
(defined in Section \ref{sec3}).

We recall that the $M$--module $\xi$ is a subquotient
of the $M$--module $V$. This means that
the object $L$ of ${\mathcal C}_X$ is a subquotient
of the object $E_V$ of ${\mathcal C}_X$. On the other
hand, $E_V$ is a strongly semistable vector bundle of
degree zero. Therefore, we conclude that
$$
\text{degree}(L)\,=\, 0\, .
$$

Thus, $E_M\, \subset\, E_{\widetilde{G}}$ is a
balanced reduction of structure group of
$E_{\widetilde{G}}$ to $M$. Now from Proposition \ref{com.}
it follows that the principal $M$--bundle $E_M$
is strongly semistable. This completes the proof
of the theorem.
\end{proof}

\begin{remark}
{\rm Assume that the fiber of the principal bundle
$E_G$, over the $k$--rational point $x$ of
$X$, admits a rational point.
If we fix a rational point in the fiber of
$E_G$ over $x$, then
$\widetilde{G}$ gets identified
with $G$, and the principal bundle
$E_{\widetilde{G}}$ gets identified with $E_G$.}
\end{remark}

If $E_H\, \subset\, E_G$ is a reduction of structure group,
to a subgroup scheme $H\, \subset\, G$, of a principal
$G$--bundle $E_G$ over $X$, then the adjoint
bundle $\text{Ad}(E_H)$ is a subgroup scheme of the
group scheme $\text{Ad}(E_G)$ over $X$.

\begin{theorem}\label{thm2}
Let $E_G$ be a strongly semistable principal $G$--bundle 
over a geometrically irreducible smooth projective
curve $X$ defined over $k$, where $G$ is a linear algebraic
group defined over $k$ with the property that $G$ does not
admit any nontrivial character which is trivial on $Z_0(G)$.
Fix a $k$--rational point $x$ of $X$.
Let $H\, \subset\, G$ be a subgroup scheme and $E_H\,
\subset\, E_G$ a balanced
reduction of structure group of $E_G$ to $H$.
Then the image in ${\rm Ad}(E_G)_{x}$
of the monodromy group scheme $M$ (image by the homomorphism
in Proposition \ref{prop.id.}) is contained in the subgroup
scheme ${\rm Ad}(E_H)_{x}\,\subset\,{\rm Ad}(E_G)_{x}$.
\end{theorem}

\begin{proof}
Take any indecomposable $H$--module $V$. Let $\chi$
be the character of $(H\bigcap Z_0(G))_0$
corresponding to the indecomposable $H$--module $V$
(see Definition \ref{de.bal.}).
Let $E_V$ denote the vector bundle over $X$ associated
to the principal $H$--bundle $E_H$ for the $H$--module
$V$. To prove the theorem it suffices to show that
$E_V$ is strongly semistable, and it is an object of the
neutral Tannakian category ${\mathcal C}_{E_G}$.
(We recall that the monodromy group $M$ is constructed from
${\mathcal C}_{E_G}$; see Definition \ref{def01}.)

Since $E_H\, \subset\, E_G$ is a balanced
reduction of structure group of $E_G$ to $H$, from
Proposition \ref{com.} we know that the principal
$H$--bundle $E_H$ is strongly semistable.
As the $H$--module $V$
is indecomposable, and the principal $H$--bundle
$E_H$ is strongly semistable, we
conclude that the associated vector bundle $E_V$ is
strongly semistable. Therefore, to complete the proof of
the theorem we need to show that $E_V$ is an object of
the neutral Tannakian category ${\mathcal C}_{E_G}$.

We recall that the group $(Z_0(G)\bigcap H)_0$
is a product of copies of the multiplicative group
${\mathbb G}_m$. Therefore, the inclusion
$$
(Z_0(G)\bigcap H)_0\, \hookrightarrow\, Z_0(G)
$$
splits \cite[\S 8.5, page 115, Corollary]{Bo}. Also,
in Lemma \ref{ch.ar.} we showed that
any character of $Z_0(H)$ arises from an indecomposable
$G$--module. Therefore, there is an indecomposable 
$G$--module $\widehat{V}$ such that
$(Z_0(G)\bigcap H)_0$
acts on $\widehat{V}$ as scalar multiplication through
the earlier defined character $\chi$ (the character
through which $(Z_0(G)\bigcap H)_0$ acts on $V$).

Let $E_{\widehat{V}}$ be the vector bundle over $X$
associated to the principal $G$--bundle $E_G$ for the
$G$--module $\widehat{V}$. Since the $G$--module
$\widehat{V}$ is indecomposable, the subgroup
$Z_0(G)$ acts on $\widehat{V}$ as scalar multiplications.
As the principal $G$--bundle $E_G$ is strongly semistable,
this implies that the associated vector bundle
$E_{\widehat{V}}$ is strongly semistable
(see Corollary \ref{cor2}).

Since both the vector bundles $E_V$ and $E_{\widehat{V}}$ are
strongly semistable, the vector bundle
\begin{equation}\label{defw}
W\, :=\, E^*_{\widehat{V}}\bigotimes E_V
\end{equation}
is also strongly semistable (see Theorem \ref{rr.t}).

As both $E_{\widehat{V}}$ and $W$ are
strongly semistable, the vector bundle
$E_{\widehat{V}}\bigotimes W$ is also
strongly semistable (see Theorem \ref{rr.t}).
Furthermore, the vector bundle $E_V$
is a subbundle of $E_{\widehat{V}}\bigotimes W
\,=\, \text{End}(E_{\widehat{V}})\bigotimes E_V$;
note that ${\mathcal O}_X$ is subbundle
of $\text{End}(E_{\widehat{V}})$, and hence
$E_V$ is a subbundle of
$\text{End}(E_{\widehat{V}})\bigotimes E_V$.
We also have
$$
\mu(E_{\widehat{V}}\bigotimes W)
\, =\, \mu (E_V)\, .
$$
Therefore, to show that the strongly semistable
vector bundle $E_V$ is an object of the Tannakian
category ${\mathcal C}_{E_G}$, it suffices to show that
the strongly semistable vector bundle $W$
is an object of the neutral Tannakian
category ${\mathcal C}_{E_G}$.

The group $(Z_0(G)\bigcap H)_0$ acts on
both $\widehat{V}$ and $V$ as multiplication by scalars
through the character $\chi$. This immediately
implies that $(Z_0(G)\bigcap H)_0$ acts trivially
on the $H$--module $\widehat{V}^*\bigotimes V$. Set
$$
H'\, :=\, H/(Z_0(G)\bigcap H)_0\, .
$$
Therefore, $\widehat{V}^*\bigotimes V$ is an $H'$--module.

Set
$$
G'\, :=\, G/(Z_0(G)\bigcap H)_0\, .
$$
Since $H'$ is a subgroup scheme of $G'$, there
is a $G'$--module $V'$ such that the $H'$--module
$\widehat{V}^*\bigotimes V$ is a subquotient
of $V'$ considered as an $H'$--module
(see \cite[page 139, Proposition 2.21(b)]{DM}).

Let
\begin{equation}\label{qu.}
V'\, \twoheadrightarrow\, Q
\end{equation}
be a quotient of the $H'$--module $V'$
such that the $H'$--module
$\widehat{V}^*\bigotimes V$ is a submodule of $Q$.

Let $E_{V'}$ be the vector bundle over $X$ associated
to the principal $G$--bundle $E_G$ for the $G$--module
$V'$. We will show that $E_{V'}$
is strongly semistable of degree zero.

For that, express the $H'$--module $V'$ as a direct
sum of indecomposable $H'$--modules. Since
$E_H\, \subset\, E_G$ is a balanced
reduction of structure group, the vector bundle over
$X$ associated to the principal $H$--bundle $E_H$
for an indecomposable $H'$--module is strongly
semistable of degree zero; see Definition
\ref{de.bal.} and Proposition \ref{com.}.
Therefore, $E_{V'}$ is isomorphic
to a direct sum of strongly semistable vector bundles
of degree zero (corresponding to a decomposition of
$V'$ as a direct sum of indecomposable $H'$--modules).
Hence $E_{V'}$ is a strongly semistable
vector bundle of degree zero. We also note that
$E_{V'}$ is an object of the Tannakian category
${\mathcal C}_{E_G}$.

Consider the $H'$--module $Q$ in eqn. \eqref{qu.}.
Let $E_Q$ denote the vector bundle over $X$ associated
to the principal $H$--bundle $E_H$ for the $H$--module
$Q$. Since the $H$--module $Q$ is a quotient of $V'$,
the vector bundle $E_Q$ is a quotient bundle of $E_{V'}$.
As the reduction $E_H\, \subset\, E_G$ is balanced,
we have $\text{degree}(E_Q)\, =\, 0$ (recall that
$Q$ is an $H'$--module). Since $E_Q$
is a quotient bundle of degree zero of the
strongly semistable vector bundle $E_{V'}$ of degree
zero, we conclude that $E_Q$ is strongly semistable.
Therefore, the vector bundle $E_Q$ is also an object
of the Tannakian category ${\mathcal C}_{E_G}$.

Finally, the vector bundle $W$ (defined in eqn. \eqref{defw})
is a subbundle of $E_Q$, because the $H'$--module
$\widehat{V}^*\bigotimes V$ is submodule of $Q$.
Since the reduction $E_H\, \subset\, E_G$ is balanced,
and $\widehat{V}^*\bigotimes V$ is an $H'$--module,
we conclude that $\text{degree}(W)\, =\, 0$.
As $E_Q$ is an object
of the Tannakian category ${\mathcal C}_{E_G}$, and
$E_Q$ is strongly semistable of degree zero, we
conclude that the subbundle $W\, \subset\, E_Q$ of
degree zero is strongly semistable, and furthermore,
$W$ is an object of ${\mathcal C}_{E_G}$.
This completes this proof of the theorem.
\end{proof}

\noindent
\textbf{Acknowledgements.} The referee gave detailed
comments how the results (of a previous version) could
be generalized. The referee also gave suggestions to
improve the exposition. We are very grateful to the
referee.

%%%%%%%%%%%%%%%%%%%%%%%%%%%%%%%%%%%%%%%%%%%%%%%%%%%%%%%%%%%%%%

\end{document}